\documentstyle[12pt,amssymb]{article}
\begin{document}

\newtheorem{thm}{Theorem}
\newtheorem{lmm}{Lemma}
\newtheorem{dfn}{Definition}

\vspace{50mm}
\title{ MEROMORPHIC TENSOR CATEGORIES}

\author{YAN SOIBELMAN}

\maketitle

\vspace{5mm}

\section{Introduction}

\subsection{} Let $A$ be an associative algebra with the unit over a
commutative ring $k$. Its ``categorical'' analog is a monoidal category
${\cal A}$, i.e. a $k$-linear category equipped with a bilinear functor
$\otimes:{\cal A}^2 \rightarrow {\cal A}$ which is associative. The
category ${\cal A}$ has also a naturally defined unit object.
Similarly, the categorical analog of a commutative associative algebra
is a braided category (in this case we have functorial isomorphisms $X
\otimes Y \rightarrow Y \otimes X$ satisfying natural properties, see
[De1]). 

Suppose that we have a $k$-module $A$ equipped with a $k$-linear
product $m:A^{\otimes 2} \rightarrow A$, $m(a \otimes b)=ab$,which is
not necessarily associative. Assume that there is an algebraic group
$G$ over $k$ acting on $A$, $a \mapsto a(g)$, which commutes with the
product: $(ab)(g)=a(g)b(g)$. Then one can try to say what is a
rational
$G$-associativity. Naturally it should be an equality like
$(a_1(g_1)a_2(g_2))a_3(g_3)=a_1(g_1)(a_2(g_2)a_3(g_3))$ valid for
Zariski open subset of $G^3$. Here we can treat both sides as elements
of the space of rational functions $Funct(G^3,A^{\otimes 3})$. One can
easily check that this ``associativity for triples'' does not imply the
associativity for products of four elements, and so on. In other words
one should consider all the spaces $Funct(G^n,A^{\otimes n})$ and
state the compatibility conditions using all of them. One can treat
rational commutativity similarly. 

This paper has arisen from the attempts to understand the categorical
analogs of these and other examples. It does not contain deep
results.It can be considered as an introduction to the axiomatics
illustrated in a few interesting examples. It contains a
 part of the notes I have been taking on the topic for more than a
year. I thank to Richard Borcherds who has convinced me to publish
these notes. I would like to mention that he has developed
independently the formalism of relaxed multi-linear categorie s and
$G$-vertex algebras (see [Bo]) which is closely related to the topic of
this paper. 

\subsection{} Another motivation for this paper came from [BD] where
the notion of pseudo-tensor category was introduced. As a special case
it gives the notions of symmetric monoidal category and operad. The
former corresponds to a representable pseudo-tensor struct ure. In this
paper we advocate a similar point of view: if we want to speak about
associativity and/or commutativity constraints depending on parameters
we should do that in the language of operations rather than objects. In
the case of a symmetric monoidal category these operations are given
by $P_I({\{ X_i \}},Y)=Hom(\otimes_I X_i,Y)$ (see [BD]). In our case we
replace sets of operations by sheaves of operations , and the
representability becomes a more sophisticated problem. The point is
that we can wor k with operations without solving it. 

Thus we naturally generalize pseudo-tensor categories in two
directions: a) introducing an operad of spaces so that operations
become sheaves (or rational or meromorphic sections of sheaves) over
these spaces ; b) introducing actions of braid groups on the operations
instead of the actions of symmetric groups $Aut(I)$ in pseudo-tensor
case. 

\subsection{}

Few typical examples include: finite-dimensional representations of
quantum affine algebras, admissible representations of $GL(n,{\bf
F})$,where ${\bf F}$ is a local field, classical chiral algebras ( see
[BD]) as well as their q-deformations.In all these
 cases one can speak about analytic (rational, meromorphic) braided (or
pseudo-braided) category. We discuss shortly these examples in the main
body of the paper. 

\subsection{Acknowledgements}

I would like to thank
A.Beilinson,R.Borcherds,P.Deligne,V.Drinfeld,G.Faltings, M.Kapranov,
V.Lyubashenko,Yu.Manin,B.Mazur,V.Schechtman for many interesting
discussions.Different parts of this work were done when I was visiting
at the IHES in Bures, Institute for Advanced Study in Princeton and the
Max-Planck Institute of Mathematics in Bonn. I thank all these
institutions for the hospitality and exellent working conditions. This
work has been partially supported by an NSF grant. 

\section{Pseudo-braided categories}

\subsection{} Suppose that for each $n \ge 1$ we are given a category
${\cal T}(n)$. Suppose that: 

a) the category ${\cal T}(1)$ has a distinguished object $e$; 

b) for any $n$, $k_1$,...,$k_n$ we are given a functor ${\cal T}(n)
\times {\cal T}(k_1) \times ... \times {\cal T}(k_n) \rightarrow {\cal
T}(k_1 +...+k_n)$ such that

$(T,T_1,...,T_n) \mapsto T(T_1,...,T_n)$. 

It is called a composition functor; 

c) the composition functors are strictly associative : 

$T(T_1(X_{11},...,X_{1k_{1}}),...,T_n(X_{n1},...,X_{nk_{n}})=(T(T_1,...,T_n))(X_{11},...,X_{nk_{n}})$; 

d) $T(e,...,e)=T$, $e(T)=T$ for any $T$. 

We will call a collection ${\cal T}=({\cal T}(n))_{n \ge 1}$ a $strict$
$monoidal$ $2-operad$.It is strict, since we use equality of functors
in c) (it can be replaced by an isomorphism of functors in non-strict
case). It is monoidal since we do not require an action of the
symmetric group $S_n$ on ${\cal T}(n)$. 

$Example$.

Let $n \ge 1$ be an integer.We denote by $\cal{T}$$(n)$ the category of
connected plane trees with $n$ external edges (=tails) equipped with
the ordered labeling from $1$ to $n$, and all the edges oriented down
to the root vertex. Recall (see [GK] or [KM ] for the details) that the
morphisms in each category are either identities or compositions of
contractions of edges. There are gluing functors between different
categories: having a tree $T\in\cal{T}$$(n)$, and a sequence of trees
$T_i \in \cal{T}$$(k_i )$ $1 \le i \le n$ one can construct a tree
$T(T_1,T_2,...,T_n)$ which is an object of $\cal{T}$$(k_1+...+k_n)$ . 

Let ${\cal T}$ be a strict monoidal 2-operad. 

\begin{dfn}. 

 A ${\cal T}$-pseudo-monoidal category is a class of objects $\cal{C}$
together with the following data: 

a) A set $P_T({\{ X_i \}},Y)$ given for any $T \in {\cal T}(n)$, a
sequence $X_1,X_2,..,X_n \in \cal{C}$ and $Y \in \cal{C}$. 

b) A map $\phi_f$: 

$P_T({\{X_i\}},Y)$$\rightarrow$$P_{T^{'}}({\{X_i\}};Y)$ given for any
morphism $f$: $T^{'}\rightarrow T$ in the category $\cal{T}$$(n)$. 
This map is functorial with respect to $f$. 

c) A composition map

$\Phi$:$P_T({\{X_i\}},Y) \times
P_{T_{1}}({\{M_{1j}\}},X_1)\times...\times
P_{T_{n}}({\{M_{nj}\}},X_n)\rightarrow
P_{T(T_1,..,T_n)}({\{M_{ij}\}},Y)$

given for any objects $T \in \cal{T}$$(n)$,$T_i \in \cal{T}$$(k_i)$
,and sequences of objects ${\{X_1,...,X_n \}}$, ${\{
M_{1j}\}},...,{\{M_{nj}\}}$,$1 \le j \le k_i$,$1 \le i \le n$.

The elements of the set $P_T$ are called $operations$, composition maps
are called compositions of operations. 

d) It is required that the composition maps are transitive with respect
to the composition in ${\cal T}$, and that for any object $X \in
\cal{C}$ there exists an element $id_X \in P_e({\{X\}},X)$ with the
natural properties of the identity map (cf.[BD] ). 

\end{dfn}

If ${\cal T}$ is a strict monoidal 2-operad from the $Example$ we will
simply say that ${\cal C}$ is a $pseudo-monoidal$ $category$. Unless we
say otherwise we consider this case. 

{\bf Remarks.}

Suppose that we have a pseudo-monoidal category $ {\cal C}$. 

a) Each category $\cal{T}$$(n)$ contains the only tree $\delta_n$
without internal edges.If the morphisms in b) are isomorphisms, one can
restate the definition above in terms of the sets $P_{\delta_{n}}$
only. 

b) Similarly to [BD] we can treat a pseudo-monoidal structure as an
extension of some categorical structure on ${\cal
C}$. Indeed, considering the set of operations corresponding to the
${\cal T}(1)$ we get morphisms in ${\cal C}$. 

\begin{dfn}. 

A pseudo-monoidal category $\cal{C}$ is called pseudo-braided if for
any element $\sigma$ in the braid group $B_n$ we have (in the notations
of the previous definition) a bijection

$\mu_{\sigma}$:$P_T({\{X_i\}};Y) \rightarrow
P_T({\{X_{\sigma(i)}\}};Y)$

where $\sigma$ acts on the set ${\{1,...,n\}}$ as the corresponding
element of the symmetric group $S_n$. 

It is required that these bijections satisfy the following properties: 

a) $\mu_{\sigma \tau}=\mu_\sigma \mu_\tau$,

$\mu_1 =id$,where $1 \in B_n$ is the unit,$\sigma, \tau$ are elements
of $B_n$. 

\vspace{2mm}

b) Compatibility with the composition maps from 2.1c), which means that
the following diagram is commutative: 

$$ {\begin{array}{ccc} P_T({\{ X_i \}}_{i=1}^n,Y) \times \prod^n_{i=1}
{P_{T_{i}}({\{ M_{ij} \}}^{k_i}_{j=1},X_i)} & \longrightarrow &
P_{T(T_{1},..,T_{n})}({\{ M_{ij} \},Y)} \\

\downarrow && \downarrow \\

 P_T({\{ X_{\sigma(i)} \}} _{i=1}^n,Y) \times
\prod^n_{i=1}P_{T_{i}}({\{ M_{i \sigma (j)} \}}^{k_i}_{j=1},X_i) &
\longrightarrow & P_{T(T_{1},...,T_{n})}({\{ M_{\sigma
(\sigma_{1},...,\sigma_{n})(ij)} \}},Y) \\ \end{array}} $$ Here $1 \le
i \le n$,$1 \le j \le k_{_{i}}$,and $\sigma(\sigma_{1},...,\sigma_{n})$
denotes the element of $B_{k_{1}+...+k_{n}}$ which is the image of the
element $\sigma \times \sigma_{1} \times... \times \sigma_{n}$ under
the natural map $B_n \times B_{ k_{1}} \times ...\times B_{k_{n}}
\rightarrow B_{k_{1}+...+k_{n}}$. Vertical arrows are actions of the
braid groups, horizontal arrows are compositions of operations. 

\vspace{3mm}

c) $\mu_\sigma$ commute with morphisms in $\cal{T}$$(n)$. 

\vspace{2mm}

d) $\mu_1$ preserves $id_X$ for any object $X$.Here $1$ is the unit in
$B_1$.

\end{dfn}
 
The notion of a functor between pseudo-monoidal or pseudo-braided
categories is defined in the natural way (see [BD] for the case of
unordered sets). We will often denote pseudo-monoidal categories simply
by ${\cal C}$. Let ${\cal C}$ be a pseudo-braided category such that
all sets $P_T$ of operations depend on the sets $I$ of vertives of $T$
only. In this case the braid groups act by permutations, and we
reproduce the notion of $pseudo-tensor$ $category$ from [BD]. We are
going to use it later in the section devoted to chiral algebras. For a
pseudo-tensor category we will denote the set of operations by $P_I$
instead of $P_T$. 

Let ${\cal O}$ be a pseudo-braided category with one object. Then we
will say that we are given a $braided$ $operad$. In particular we can
speak about functors from a braided operad to a pseudo-braided
category. 

\subsection{}

We are going to generalize the notion of ${\cal T}$-pseudo-monoidal
category in the following way. We assume that for each object $T \in
{\cal T}(n)$ we are given a ringed space $(S_T,O_{S_T})$ (topological
space equipped with a sheaf of commutative rin gs), and this
correspondence is functorial for each category ${\cal T}(n)$ and
compatible with the composition functors. In particular we have a
morphism of spaces $S_T \times \prod_{i}{S_{T_{i}}} \rightarrow
S_{T(T_{1},...,T_{n})}$ corresponding to the composition of objects
$T_i$ and $T$. Thus we get a family of spaces ${\cal S}=(S_T)_{T \in
{\cal T}(n)}$, $n \ge 1$. We say that ${\cal S}$ is a ${\cal
T}$-operad.. 

Suppose that ${\cal T}$ is a strict 2-operad of plane trees, and ${\cal
S}$ is a ${\cal T}$-operad. 

\begin{dfn}. 

a) A pseudo-monoidal category over ${\cal S}$ consists of the following
data: 

(i) a class of objects ${\cal C}$; 
 
(ii) for each plane tree $T \in \cal {T}$$(n)$, a sequence of objects
${\{ X_i \}}_{i=1}^{i=n}$ and an object $Y$ we are given a sheaf (of
sets, vector spaces, etc) $P_T( S_T;{\{ X_i \}},Y)$ such that the axioms
of the Definition 1 are satisfied (with the obvious change ``sets'' to
``sheaves''; here $id_X$ must be a section of the sheaf $P_e(S_e;{\{X
\}},X)$). 

\vspace{2mm}

b) A pseudo-monoidal category over ${\cal S}$ is called pseudo-braided
if for each $ \cal {T}$$(n)$ we have an action of the braid group $B_n$
on the sheaves of operations such that the properties of the Definition
2 are satisfied ($B_n$ does not act on t he base space $S_T$). 
\end{dfn}

The notion of a functor between pseudo-monoidal or pseudo-braided
categories over ${\cal S}$ is defined in the natural way. If we have
two pseudo-monoidal categories defined over different operads of
spaces, then a functor from one to another consists of a morphism of
the ${\cal T}$-operads of spaces and morphisms of the sheaves of
operations compatible with it. We leave to the interested reader to
write down all the diagrams. 

\vspace{4mm}

{\bf Remarks}. 

a) If all the spaces are points we recover the previous definitions.
Sometimes we will omit ${\cal S}$ simply saying that we have a
pseudo-monoidal or pseudo-braided category. We hope it will not be
confusing.

b) If $\cal {F}$ is a sheaf over $X$ , $\cal {G}$ is a sheaf over $Y$,
and $f:X \rightarrow Y$ is a morphism of spaces then we say that a
morphism of the sheaves is a morphism $\cal {F}$$ \rightarrow
f^{*}$$\cal {G}$. In particular, since our spaces are ringed, the
composition maps are $\cal {O}$$_{{S_{T(T_{1},...T_{n})}}}$-linear (cf. 
[GK] for the notion of a sheaf on a topological operad).  \vspace{3mm}

\subsection{}

Let ${\cal C}$ be a pseudo-monoidal category over ${\cal S}$.  Assume
that our spaces are reduced irreducible schemes or reduced irreducible
complex analytic spaces. Then for each space $S_T$ the sheaf of
operations is a sheaf of quasi-coherent $\cal {O}$$_{S_{T}}$-modules
where  ${\cal O}_X$ denotes the structure she af of $X$. We can extend
the space of operations on $S_T$ to obtain a vector space over the
field of rational functions on $S_T$ (for schemes) or over the field of
meromorphic functions on $S_T$ (for complex analytic spaces). Let us
denote either of these
 fields by $K(S_T)$. Clearly these vector spaces define a new
pseudo-monoidal structure on ${\cal C}$. We will call it $localization$
of the original one. The pseudo-braided case can be treated similarly. 

Suppose now that we are given an operad of spaces as above, the vector
spaces $P_T(S_T;{\{ X_i \}},Y)$ over $K(S_T)$ as well as maps $\phi_f$
and composition maps between them (see Definition 1). 

\begin{dfn}. 

a)If the conditions of the Definition 1 are satisfied we say that
${\cal C}$ is: 

a rational pseudo-monoidal category if ${\cal S}$ is an operad of
schemes; 

a meromorphic pseudo-monoidal category if ${\cal S}$ is an operad of
complex analytic spaces. 

b)If in addition we have actions of braid groups as in the Definition 2
we say that ${\cal C}$ is a rational (resp.meromorphic) pseudo-braided
category. 

\end{dfn}

\vspace{2mm}

{\bf Remark}. 

Intuitively one can imagine a pseudo-monoidal category over ${\cal S}$
as given by a class of objects and spaces of operations which are
families parametrized by the spaces $S_T$. This parametrization is
either complex analytic or algebraic regular depen ding on the category
of spaces. Meromorphic or rational structures are given by similar
families which might have singularities . In the pseudo-braided case
the action of the braid group can be singular itself. 

 \subsection{Example}

Let ${\cal C}$ be a {\bf C}-linear category, $G$ be a complex analytic
group acting on the objects from the left:$X \mapsto X(g)$ for any $X
\in {\cal C}$ and $g \in G$. 

We can define an operad of spaces in the following way. To each plane
tree with $n$ tails we assign the group $G^n$, $n \ge 0$, $G^0$ is the
trivial group. The compositions are given by the maps $G^n \times
G^{k_{1}} \times ... \times G^{k_{n}} \rightarrow G^{k_{1}+...+k_{n}}$
such that
$(g_1,...,g_n;g_{11},...,g_{1k_{1}};...,g_{n1},...,g_{nk_{n}}) \mapsto
(g_{1k_{1}},...,g_{nk_{n}})$.  Suppose that for each $G^n$ we have
families of complex vector spaces $P_T(G^n;{\{ X(g_{i}) \}},Y)$
parametrized by its points. Assume that in such a way we obtain
analytic sheaves over $G^n$ satisfying the condition a) (resp.b)) of
the Definition 4. Thus we get a pseudo-monoidal (resp. pseudo-braided)
category over ${\cal
 G}=(G^n)_{n \ge 0}$. If ${\cal C}$ is a pseudo-braided category we get
a pseudo-braided category over ${\cal G}$. We call it
pseudo-monoidal(resp.pseudo-braided) $G-category$ if these sheaves are
equivariant with respect to the left action of $G$ on $G^ n$ (in
pseudo-braided case the action of the braid group also has to be
compatible with this equivariance).  If the conditions of the
Definition 5 are satisfied then we obtain a $meromorphic$
$pseudo-monoidal$ $G-category$ (or its pseudo-braided version). This
construction can be done in a pure algebraic setting as well ($G$ has
to be an algebraic group in this
 case). We will see that many examples which naturally appear in
practice can be obtained in this way. 

\subsection{} In this subsection we discuss the case of representable
pseudo-monoidal and pseudo-braided structures. We restrict ourselves
mostly to the case of pseudo-monoidal categories.The pseudo-braided
case is similar. 

Let us make a few comments. If we have a usual pseudo-monoidal category
${\cal C}$ then we say that its pseudo-monoidal structure is
representable if for any plane tree $T \in {\cal T}(n)$ there exists a
functor $F_T: {\cal C}^n \rightarrow {\cal C}$ such
 that we have a functorial isomorphism of sets (vector
spaces, modules,etc.) $P_T({\{ X_i \}},Y) \rightarrow Hom(F_T({\{ X_i
\}}),Y)$ valid for any sequence of objects $X_1,X_2,...,X_n ,Y$ of
${\cal C}$. It is required to be compatible with the morphisms i n
${\cal T}(n)$ and with the gluing of trees and composition of functors
in $(Funct({\cal C}^n,{\cal C}))_{n \ge 1}$. 

A typical example is given by a monoidal category.Then each binary
plane tree $T$ gives rise to a functor (tensor product with the maximal
bracketing prescribed by $T$). The Maclane coherence theorem allows us
to extend this construction to get a functor $F_T$ for each plane tree
$T$. Similarly a usual braided category gives an example of a
representable pseudo-braided structure. 

Suppose now that we have a pseudo-monoidal category ${\cal C}$ over
${\cal S}$. Then to say that it is representable we need to speak about
families of objects of ${\cal C}$ parametrized by the spaces $S_T$. We
briefly recall an appropriate language of sheaves of categories. 

Let $X$ be a topological space (or a Grothendieck topology). We recall
(see [Gr] for the details) that a sheaf (=stack)  ${\cal F}$ of
categories on $X$ is given by: 

a) a presheaf of categories.It assigns a category to each open set from
the category of open sets on $X$ in such a way that the usual
properties hold; 

b) gluing (=descent) data for ${\cal F}$ making this presheaf into a
sheaf. 
 
In particular for each open set $U$ we have a fiber category ${\cal
F}_U$, and for each morphism of open sets $f:V \rightarrow U$ there is
a pull-back functor $f^{*}:{\cal F}_U \rightarrow {\cal F}_V$. These
structures are functorial in the natural way. If all fiber categories
are subcategories of a category ${\cal C}$ we say that we have a sheaf
on $X$ with values in ${\cal C}$. If all fiber categories coinside with
${\cal C}$ we have a constant sheaf with the fiber ${\cal C}$. 

 A sheaf of categories is called $locally$ $constant$ if it is locally
isomorphic to a constant one.It is called also a $local$ $ system$ with
values in ${\cal C}$. 

In the rest of this subsection we will assume that all $S_T$ are
complex analytic spaces unless we say otherwise. The case of schemes
can be treated similarly. 

Let ${\cal C}$ be a ${\bf C}$-linear category. Then for any complex
analytic space $X$ we can construct a free sheaf of categories over
$X$,namely ${\cal O}_X \otimes {\cal C}$. We say that a sheaf ${\cal
F}$ of categories over an analytic space ${\cal M }$ is a $bundle$ with
the fiber ${\cal C}$ if for any open $U \subset {\cal M}$ there is a
morphism of the open subspaces $f: V \rightarrow U$ in ${\cal M}$ such
that the the pull-back category $f^{*}({\cal F}_U)$ is isomorphic to
${\cal O}_V \otimes {\cal C}$. In particular any object $A \in {\cal
C}$ can be viewed as a category with one object and hence defines a
bundle with the fiber $A$ : to each open $U \subset {\cal M}$ it
assigns ${\cal O}(U) \otimes A$. There is a natural functor from the
category of local systems with a ${\bf C}$-linear fiber ${\cal C}$ to
the category of bundles with the fiber ${\cal C}$. 

Let ${\cal C}$ be as above, ${\cal S}$ be an operad of analytic spaces.
Suppose that for each tree $T \in {\cal T}(n)$ we are given a functor
$F_T \in Funct({\cal C}^n,{\cal C})$ such that the correspondence $T
\mapsto F_T$ is compatible with the morphisms and composition of trees
(one can say that we are given a morphism of strict monoidal 2-operads
${\cal T} \mapsto {\{ Funct( {\cal C}^n,{\cal C}) \}}_{n \ge 1}$). Then
for each $T$ and sequence of objects ${\{X_i \}}_{1 \le i \le n}$ we
have a trivial
 bundle over $S_T$ isomorphic to ${\cal O}_{S_T} \otimes Hom(F_T({\{
X_i \}}),Y)$ and these bundles form a pseudo-monoidal operad over
${\cal S}$.We denote each bundle by ${\cal F}_T$.It is a bundle of
categories with the fiber $F_T \in Funct({\cal C}^n,{ \cal C})$. 

Assume now that we are given a pseudo-monoidal structure on ${\cal C}$
extending its categorical structure (see Remark in 2.1). 
  
\begin{dfn}. 

We say that the pseudo-monoidal structure of ${\cal C}$ is
representable if for each tree $T \in {\cal T}(n), n \ge 1$,it is
locally isomorphic to the bundle of categories ${\cal F}_T$ on ${\cal
S}_T$ in such a way that we have a morphism of pseudo-monoid al
categories over ${\cal S}$.This means that the following a) and b)
hold: 

a) $P_T(S_T;{\{ X_i \}},Y)$ is locally isomorphic (as an analytic
sheaf) to ${\cal O}_{S_T} \otimes Hom(F_T({\{ X_i \}}),Y)$

(in self-explaining notations); 

b) these isomorphisms are compatible with morphisms in ${\cal T}(n)$,
with the operadic structure on ${\cal T}$ ,composition maps between the
sheaves of operations and composition of functors $F_T$. 

c) Suppose that all spaces $S_T$ are reduced and irreducible. We say
that we have a representable meromorphic pseudo-monoidal structure if
a) and b) hold for the corresponding spaces over the fields of
meromorphic functions. 

(In particular $P_T(S_T;{\{X_i \}},Y)$ is isomorphic to $K(S_T) \otimes
Hom(F_T({\{ X_i \}}),Y)$). 

d) Suppose in addition to c) that all morphisms ${\phi_f}$ (see
Definition 1) and all composition maps are isomorphisms. Then we say
that our category is rational monoidal (in the case of schemes) or
meromorphic monoidal (in the case of complex analytic spaces). In this
case the functor $F_{\delta_2}$ , ${\delta_2} \in {\cal T}(2)$ (see
Remark in 2.1 about $\delta_n$) is called a tensor product. In case if
the representable structure was pseudo-braided we would call ${\cal C}$
a rational (or meromorphic) 
 braided category.

\end{dfn}

Intuitively we can think about a meromorphic braided category as about
a class of objects ${\cal C}$ such that if $X$ and $Y$ are two of them
then $Hom(X,Y)$ is a vector space over the field of meromorphic
functions on $S_{\delta_1}$. The tensor product o f $X$ and $Y$ is a
vector space $X \otimes Y$ over the field of meromorphic functions on
$S_{\delta_2}$. The generator of the braid group $B_2$ gives an
isomorphism of this vector space with the similar vector space $Y
\otimes X$. The same can be done for any tree $T \in {\cal T}(n)$ (in
this case we have an action of $B_n$ of course).Compatibility with the
gluing of trees means that our tensor product gives rise to an
``associativity constraint with parameters'' which might fail to be an
isomorphism on
 some closed subspaces of $S_T$. Similarly the action of $B_n$ gives
rise to a ``commutativity constraint with parameters'' which can fail
to be an isomorphism on some (other) closed subspaces of $S_T$.
\vspace{2mm}

{\bf Remarks}. 

a) This definition can be generalized to a quite general framework of
Grothendieck topologies. 

b) In the braided case we have an action of $B_n$ on $Funct({\cal
C}^n,{\cal C})$ through the natural homomorphism $B_n \rightarrow S_n$(
the symmetric group $S_n$ acts on the functors via permutations of
variables). This action is required to be compatible with the
corresponding action on the operations of a given pseudo-braided
category. It is clear how to generalize this picture to the meromorphic
case.We are going to use all these generalizations without extra
comments.

\subsection{The case of braided categories}

Let ${\cal C}$ be a braided monoidal category (see for ex.[De 1]). Then
to any binary plane tree $T \in {\cal T}(n)$ we can assign a functor
$F_T: {\cal C}^n \rightarrow {\cal C}$ which maps a sequence of objects
${\{ X_i \}}_{i=1}^{i=n}$ to their tensor
 product with the bracketing prescribed by $T$. 

Let ${\overline{M_{0,n+1}}}$ be the moduli space of complex stable
curves of genus $0$ with $n+1$ marked points and a non-zero tangent
vector assigned to the last point. The real strata of this moduli space
are parametrized by the elements of ${\cal T}(n ) \times S_n$. Binary
trees correspond to the zero-dimensional strata and the tree $\delta_n$
corresponds to the open strata. Following [De 2] one can construct a
local system of categories on ${\overline{M_{0,n+1}}}$ with values in
$Funct({\cal C}^n,{\cal C})$ . To do this one uses the associativity
constraint in ${\cal C}$. For example let $f:T \rightarrow \delta_3$
and $g:T^{\prime} \rightarrow \delta_3$ are the only non-trivial
morphisms in ${\cal T}(3)$. For the binary trees $T$ and $T^{\prime}$
we have the corresponding functors $F_T$ and $F_{T^{\prime}}$ (tensor
products of three objects with two possible bracketings). Then we have
an isomorphism of functors $F_T \rightarrow F_{T^{\prime}}$ given by
the associativity constraint. From the geometric
 viewpoint we have two embeddings of the real strata corresponding to
$(T, 1)$ and $(T^{\prime},1)$ to the boundary of ${\overline{M_{0,4}}}$
(here $1$ is the unit
 of the symmetric group). One can construct a local system on the open
stratum corresponding to $(\delta_3,1)$ in such a way that it has a
constant fiber $F_T$ and its specialization to the stratum
corresponding to $(T^{\prime},1)$ is identified with $F_{ T^{\prime}}$
via the associativity constraint. In general, to get a local system on
the moduli space ${\overline{M_{0,n+1}}}$ we use the action of the
braid group (remark that we have $n!$ real components for
${\overline{M_{0,n+1}}}$ corresponding to a fixed tree). The braid
group $B_n$ acts on $F_T$, $T \in {\cal T}(n)$ via the commutativity
constraint. 
 On the other hand it is a fundamental group of the corresponding
moduli space. Note that our local systems are compatible with the
morphisms of trees. Namely if $f:T^{\prime} \rightarrow T$ is a
morphism in ${\cal T}(n)$ then we have a closed embedding o f the
component of the moduli space ${\overline{M_{0,n+1}}}$ corresponding to
$T$ to the boundary of the component corresponding to $T^{\prime}$.
Specialization of the local system to the boundary component is
well-defined. This specialization is a local system isomorphic to the
existing one (as above: the action of the braid group is interpreted as
an action of the fundamental group ). Gluing of trees is compatible
with the embeddings of the products of smaller components to the larger
ones. In the
 language of tensor products gluing of trees corresponds to the
composition of functors. 

According to [De 2] there is one-to-one correspondence between braided
categories and local systems on the operad of spaces
$({\overline{M_{0,n+1}}})_{n \ge 1}$ equipped with compatibilities
described above. The corresponding vector bundles of $Hom$'s give an
example of a representable pseudo-braided structure over
$({\overline{M_{0,n+1}}})_{n \ge 1}$. Of course this pseudo-braided
category is equivalent to a usual pseudo-braided category (over a
point) which is also representable. The equivalence is given
essentially by the restriction of the local systems to the
zero-dimensional strata. This observation explains why braided
categories often appear in the form of families of local systems.  A
famous example is the family of conformal blocks of the WZW model in
Conformal Field Theory (see [MS]). They appear as the local systems of
solutions of Knizhnik-Zamolodchikov equations and can be thought as
local systems on ${\overline{M_{0,n+1}}}$. These local systems give
rise to a pseudo-braided category over $({\overline{M_{0,n+1}}})_{n \ge
1}$. On the other hand, according to [Dr] and [KL] this pseudo-braided
category is representable, and moreover equivalent to a usual braided
category (in [KL] it is a category of certain highest weight
representations of a simply-laced affine Kac-Moody algebra with the
fixed central charge. The tensor product is the so-called ``fusion''
tensor product. It corresponds to the operator product expansion in
Conformal Field Theory).

\subsection{ Remark on meromorphic braided $G$-categories}

Suppose that we are in the assumptions of the Example 2.4, and moreover
we have a meromorphic monoidal category over ${\cal G}$. We say that we
are given a $meromorphic$ $monoidal$ $G-category$ if the action of $G$
commutes with the tensor product :$(X \otimes Y)(g)=X(g) \otimes
Y(g)$. We define meromorphic $G$-braided categories in a similar
fashion. 

A typical situation when meromorphic monoidal $G$-categories can appear
is the following. Suppose that ${\cal C}$ is a ${\bf C}$-linear
category with the tensor product i.e. with a functor $\otimes: {\cal
C}^2 \rightarrow {\cal C}$. It can happen that our tensor product is
not necessarily associative or commutative. Assume that there is an
analytic group $G$ acting on the category ${\cal C}$ in such a way that
it commutes with the tensor product. It can happen that for any objects
$X$,$Y$,$Z$ we have functorial isomorphisms
$a_{X(g_1),Y(g_2),Z(g_3)}:(X(g_1) \otimes Y(g_2)) \otimes Z(g_3)
\rightarrow X(g_1) \otimes (Y(g_2) \otimes Z(g_3))$ which is
meromorphic on $G^3$ as well as similar isomorphisms for higher
iterations of the tensor product. If all these isomorphisms are
compatible with the compositions of the tensor product functors we
say that we have a meromorphic monoidal $G$-category. We remark that
the ``higher'' associativity constraints for the iterated tensor
products should be given as a part of the data. This differs from the
case of the usual monoidal categories. 

Similarly, it can happen that we do not have an isomorphism $X \otimes
Y$ and $Y \otimes X$ but we have a meromorphic isomorphism
$c_{X(g_1),Y(g_2)}:X(g_1) \otimes Y(g_2) \rightarrow Y(g_2) \otimes
X(g_1)$ which is compatible with the meromorphic associativity. Then
we have a meromorphic $G$-braided structure on ${\cal C}$. We will
consider an interesting example in the next section. 

\section{Quantum affine algebras}

\subsection{} Let $g$ be a complex finite-dimensional simple Lie
algebra. We fix an invariant bilinear form on it. Then for a given
non-zero complex number $q$ ($q$ is not a root of $1$) one can define
the Drinfeld-Jimbo quantized universal enveloping algebra $U_q(g)$
which is a complex Hopf algebra. Let $g^{\prime}$ be the affine
non-twisted Kac-Moody algebra corresponding to $g$. The corresponding
quantized universal enveloping algebra $U_q(g^{\prime})$ is a Hopf
algebra containing $U_q(g)$. We are not going to repro duce its
description here reffering the reader to [L] and [KS].Following [L] we
will denote the generators of $U_q(g^{\prime})$ by $E_i$ ,$F_i$ ,
$K_{\mu}$. Here $\mu$ belongs to the co-weight lattice $\Lambda^{\vee}$
of $g^{\prime}$,$i$ runs through the finite set $I$ corresponding to
the affine irreducible root datum of $g^{\prime}$.We denote by $i_{0}$
the special vertex of $I$. Then the Hopf subalgebra $U_q(g)$ is
generated by the subset of the above-mentioned generators for which $i$
runs through $I -{\{ i_{0} \}}$ , and $\mu$ belongs to the co-weight
lattice of the Lie algebra $g$. 

We recall that there is a natural group homomorphism ${\bf
Z}$$[I]$$\rightarrow \Lambda$ where $\Lambda$ is the weight lattice of
$g^{\prime}$. The kernel is generated by the element $\sum_{i \in
I}{n_i i}$ such that $n_{i_{0}}=1$.We denote by $K_i$ the g enerator
$K_{\mu}$ with $\mu$ equals to $(i \cdot i)i/2$. Then the element
$Z$=$\prod_{i \in I}{K_{i}^{n_{i}}}$ is central in $U_q(g^{\prime})$. 
For simplicity we will assume that $g$ is simply-laced. 

\subsection{} We denote by $\cal{C}$ the category of finite-dimensional
unital $U_q(g^{\prime})$-modules such that $Z$ acts on the objects as
the identity morphism and all $K_{\mu}$ are diagonalizable with the
eigen-values from $q^{\bf Z}$. The category $\cal{C}$ is a monoidal
category with the associativity constraint to be identity. There is an
action of the group ${\bf C^{*}}$ on the category $\cal{C}$. It comes
from the corresponding action on $U_q(g^{\prime})$ such that $E_i
\rightarrow zE_i$, $F_i\rightarrow z^{- 1}F_i$, and all $K_{\mu}$ are
stable under the action. 

This action makes an object $X$ into an object $X(z)$.It is clearly
compatible with the tensor product in ${\cal{C}}$. It was shown in [KS]
that the universal quantum $R$-matrix of $U_q(g^{\prime})$ defines a
family of morphisms $c_{X(z_{1}),Y(z_{2})}:X(z _1) \otimes Y(z_2)
\rightarrow Y(z_2) \otimes X(z_1)$ which is meromorphic in $z_1/z_2$. 
The following result was proved in [KS]. 

\vspace{2mm}

\begin{thm}. 

The category ${\cal{C}}$ carries a structure of a meromorphic braided
${\bf C^{*}}$-category. For any two objects $X$ and $Y$ the square of
the commutativity constraint

$c_{X(z_1),Y(z_2)}c_{Y(z_2),X(z_1)}$ is an elliptic function on the
curve $\cal{E}$=${\bf C^{*}}/q^{\bf 2Z}$ with values in the vector
space $End_{\bf C}(X \otimes Y)$. 

\end{thm}

\subsection{}

In this subsection we are going to use the categories ${\cal O}_z^{+}$
of smooth $U_q(g^{\prime})$-modules with the fixed central charge $z$
defined in [KS], Section 3. Weyl modules (induced from
finite-dimensional simple $U_q(g)$-modules ) are examples of the
objects of ${\cal O}_z^{+}$. Let $End_p$ be the category of
endomorphisms of ${\cal O}_z^{+}$ for $p=zq^h$ where $h$ is the dual
Coxeter number of $g$. We assume that $\vert p \vert $ is greater than
1. Obviously $End_p$ is a monoidal category (tens or product is given
by the composition of functors). We can treat it as a meromorphic
monoidal $G$-category for any analytic group $G$. We choose $G$=${\bf
C^{*}}$. 

Then one of the main results of [KS] can be reformulated as follows. 

\begin{thm}. 

There is a functor $F:{\cal C} \rightarrow End_p$ of meromorphic
monoidal ${\bf C^{*}}$-categories. 

Define $F(z)$ as $F(X(z))$ for any $X \in {\cal C}$, $z \in {\bf
C^{*}}$. Then the functor $F(z)$ is isomorphic to $F(zp^2)$. Therefore
the family ${\{ F(z) \}}$ gives rise to a ``line bundle'' over the
elliptic curve ${\bf C^{*}}/p^{\bf 2Z}$.

\end{thm}

\subsection{Remarks on the Yangian case}

For any simple complex Lie algebra $g$ V.Drinfeld has constructed an
infinite-dimensional Hopf algebra $Y(g)$ called the Yangian of $g$ (see
[Drinfeld,ICM-86]). The Hopf algebra $U_q(g^{\prime})$ can be
considered as a quantization of a centrally extend loop algebra
$g[t,t^{-1}]$. Similarly the Yangian $Y(g)$ can be considered as a
quantization of the regular loop algebra $g[t]$. Let us denote by
$I_{\mu}$ and $J_{\mu}$ the standard generators of $Y(g)$ (see [D1]).
Then there is an action of the additiv e group ${\bf C}$ on the Hopf
algebra $Y(g)$ such that $I_{\mu} \mapsto I_{\mu}$ and $J_{\mu} \mapsto
J_{\mu}+zI_{\mu}$ for any $z \in {\bf C}$.Therefore ${\bf C}$ acts on
$Y(g)$-modules. 

Theories of finite-dimensional modules of $U_q(g^{\prime})$ and $Y(g)$
are completely parallel (see [D2] ). In particular one has the
following result. 

\begin{thm} The category of finite-dimensional $Y(g)$-modules is a
meromorphic braided ${\bf C}$-category. The square of the commutativity
constraint (cf. Theorem 1) is a ${\bf Z}$-periodic meromorphic function
on ${\bf C}$. 

\end{thm}

One can also obtain an analog of the Theorem 2. For this one should
consider a central extension of the double of the Yangian $D(Y(g))$
(see [S]). 

\section{Representations of $GL({\bf F})$ and Eisenstein series}

This section is based on the unpublished manuscript by M.Kapranov. 
\subsection{} Let ${\bf F}$ be a local field.We are going to consider
admissible representations of the groups $GL(n,{\bf F})$,$n \ge 0$ (see
[BZ]). 

We recall that if $V_i$,$i=1,2$ are admissible representations of
$GL(n_i,{\bf F})$ then one can define a new admissible represntation of
$GL(n_1+n_2,{\bf F})$ by the formula

$$ V_1 \odot V_2=Ind^{GL(n_1+n_2,{\bf F})}_{P(n_1,n_2)}(V_1 \otimes
V_2)  $$

where

$$ P(n_1,n_2)=\left( \begin{array}{cc} GL(n_1,{\bf F}) & \ast \\ 0 &
GL(n_2,{\bf F}) \\ \end{array} \right).  $$

\begin{dfn}. 

 $GL({\bf F})$- module is a collection $V=(V_n)_{n \ge 0}$ of
admissible representations such that $V_n$ is a representation of
$GL(n,{\bf F})$.  \end{dfn}

Proof of the following proposition is straightforward. 

\begin{thm} The following operation makes the category ${\cal B}$ of
$GL({\bf F})$-modules into a monoidal category $$ V \odot W= \oplus_n
(V \odot W)_n $$ where $$ (V \odot W)_n=\oplus_{i+j=n}V_i \odot W_j $$
\end{thm}

Let us remark that there is an action of the group ${\bf C}$ on the
$GL(n,{\bf F})$-modules , $V \mapsto V(z)$ where

$$ V(z)=V \otimes |det|^z $$

This action gives rise to the action of the group ${\bf C^{\infty}}$ on
the category ${\cal B}$ where ${\bf C^{\infty}}$ is the infinite
product of the additive groups ${\bf C}$.  We recall (see loc.cit) that
if $V_i$ are admissible representations of $GL(n_i,{\bf F})$,$i=1,2$,
then for generic complex numbers $z_i$, $i=1,2$ there is an
intertwining operator:  $$ A_{V_1(z_1),V_2(z_2)}: V_1(z_1) \odot
V_2(z_2) \rightarrow V_2(z_2+n_1) \odot V_1(z_1-n_2)  $$ Let us define
a new tensor product of the above-mentioned representations:  $$ V_1
\bullet V_2=V_1(n_2/2) \odot V_2(-n_1/2)  $$ Then it is easy to see
that the operator
$M_{V_1(z_1),V_2(z_2)}=A_{V_1(z_1+n_2/2),V_2(z_2-n_1/2)}$ defines an
intertwiner $V_1(z_1) \bullet V_2(z_2) \rightarrow V_2(z_2) \bullet
V_1(z_1)$ for generic $z_1,z_2$.  Using the same definition as for
$\odot$ we extend the tensor product $\bullet$ to the category ${\cal
B}$.  \begin{thm} This makes ${\cal B}$ into a meromorphic braided
${\bf C^{\infty}}$-category.  \end{thm}

Here we understand meromorphic function as being such when restricted
to any finite product of analytic groups ${\bf C}$. 

Now we need to recall some facts about Yangians and affine quantum
algebras. 

Let us fix $l$ and $n$ and consider all $Y(sl(n))$-modules $M$ whose
$sl(n)$-irreducible components appear in $({\bf C^n})^{\otimes l}$.
Suppose that $l$ is the minimal number with this property. In this case
we say that $M$ has level $l$. Every object in the category of
$Y(sl(n))$-modules has some level. Let us consider the category ${\cal
A}({n})$ formed by the sequences $(M_i)_{i \ge 0}$ of
$Y(sl(n))$-modules, such that $M_0=0$, each $M_i$ has level less or
equal than $i$. 

 Replacing $Y(sl(n))$ by $U_q(sl(n)^{\prime})$ where $sl(n)^{\prime}$
is the non-twisted affine Lie algebra corresponding $sl(n)$ we obtain
the category ${\cal E}(n)$ instead of ${\cal A}(n)$. In this case we
use the tensor power of the natural representa tion of $U_q(sl(n))$ to
define the level.  We refer the reader to [CP] for the details. It can
be shown that ${\cal A}(n)$ is a meromorphic braided ${\bf
C^{\infty}}$-category, and ${\cal E}(n)$ is a meromorphic braided
$({\bf C^*)^{\infty}}$-category. 

Similarly to the category ${\cal B}$ one can define the categories
${\cal H}$ and ${\cal L}$. The first one consists of the sequences
$(V_i)_{i \ge 0}$ such that $V_0=0$, and each $V_i$ is a module over
the affine Hecke algebra $H_i$. In the second case we use sequences of
modules over degenerate Hecke algebras $\Lambda_i$ defined in [D3]. 
Both categories can be equipped with tensor products (each tensor
product is similar to the parabolic induction in the $GL({\bf
F})$-case. It is called the $Zelevins ky$ $tensor$ $ product$ in [CP]). 

The following result can be derived from [D3],[CP]. 

\begin{thm} a) The category ${\cal H}$ is a meromorphic braided $({\bf
C^*})^{\infty}$-category. The category ${\cal L}$ is a meromorphic
braided ${\bf C}^{\infty}$-category. 

b) There is an equivalence of ${\cal H}$ to the subcategory of ${\cal
E}(n)$ consisiting of those sequences $(V_i)$ for which $V_i$ has
exactly level $i$ for $i<n$. This equivalence is compatible with the
structures of meromorphic braided $({\bf C^*})^{\ infty}$-categories. 

Similar result holds for the pair ${\cal L}$ and ${\cal A}(n)$. 

\end{thm}

 \subsection{} The results of the previous subsection can be naturally
extended to the global case. Namely let $k$ be a global field, ${\bf
A}$ be its ring of adeles. A representation of $GL(n,{\bf A})$ is
called $automorphic$ if it can be embedded into the regular
 representation in the space C(GL(n,{\bf A})/GL(n,k)).Let $M=\otimes_{x \in
X} M_x$ be a $GL(n,{\bf A})$-module (here $X$ denotes the set of places
of $k$). Using the adelic norm we can ``twist'' $M$ by $|det|^s$. This
gives an action of the additive group. Le t ${\cal G}$ be the the
category consisting of sequences $(M_i)_{i \ge 0}$ such that $M_i$ is a
$GL(i,{\bf A})$-module, ${\cal J}$ be the subcategory of automorphic
modules. Using the ``local'' definitions from the subsection 1 we
define tensor products $\odot$ and $\bullet$ on ${\cal G}$. It is known
that ${\cal J}$ is closed under these tensor products. Using the
twisting by $|det|^s$ (here $|a|$ denotes the adelic norm of $a$) we
make ${\cal G}$ and ${\cal J}$ into ${\bf C}^{\infty}$-categories. Let
u s fix $\bullet$ as the tensor product on ${\cal G}$. 

\begin{thm} a) The category ${\cal G}$ becomes a meromorphic braided
${\bf C}^{\infty}$-category. 

b) The subcategory ${\cal J}$ is a meromorphic tensor ${\bf
C}^{\infty}$-subcategory (i.e. the square of the commutativity
constraint is $1$). 

\end{thm}

The Eisenstein series construction defines a linear map $C(GL(n,{\bf
A})/GL(n,k)) \bullet C(GL(m,{\bf A})/GL(m,k)) \rightarrow C(GL(n+m,{\bf
A})/GL(n+m,k))$,

$$ Eis(f)=\sum_{\gamma \in GL(n+m,k)/P(n,m)} {f(g \gamma)} $$

Let us consider ${\cal L}=\{ C(GL(n,{\bf A})/GL(n,k)) \}_{n \ge 0}$ as
an object of ${\cal J}$. Then we have constructed a morphism $Eis$:
${\cal L} \bullet {\cal L} \rightarrow {\cal L}$. 

\begin{thm} This makes ${\cal L}$ into an associative algebra in the
meromorphic tensor category ${\cal J}$.  \end{thm} \vspace{2mm}

{\bf Remarks}

a) We recall that we can consider a functor from any operad to a
meromorphic pseudo-braided category. Taking for example the associative
operad $Ass$ we arrive (cf. [BD]) to the notion of an associative
algebra in a meromorphic tensor category. 

b) We can restrict ourselves to the subgroup ${\bf C}$ diagonally
embedded into ${\bf C}^{\infty}$. Then we obtain meromorphic functions
in one variable. In particular part b) of Theorem 7 leads to
functional equations for $L$-functions. 

c) One can establish equivalences of the meromorphic tensor categories
from Theorems 5 and 6 and a meromorphic tensor subcategory of ${\cal
G}$ corresponding to the weakly ramified representations. We leave to
the reader this reformulation of the well-known results connecting
representations of Hecke algebras and groups $GL$ over local fields.

\section{Classical and quantum chiral algebras}

\subsection{}

In this subsection we follow [BD]. We are going to use a special case
of the notion of pseudo-braided category,namely a pseudo-tensor
category. We recall that the sets of operations for a pseudo-tensor
category depend on the unordered sets of vertices o f the trees. Let
$X$ be a smooth complex curve. For any finite $I$ we denote by
$\Delta^{(I)}$ the diagonal embedding $X \hookrightarrow X^I$. Let
$j^{(I)}:U^{(I)} \hookrightarrow X^I$ be the embedding of the
complement of the diagonal divisor. Let ${\cal
 M}(X)$ be the category of right ${\cal D}$-modules on $X$. Then for a
sequence of objects ${\{ L_i \}}_{ i \in I}, N$ of this category one
can define the following sets

$$ a)P_I^{*}({\{ L_i \}},N)=Hom(\boxtimes L_i,\Delta_{*}^{(I)}N);  $$

$$ b)P_I^{ch}({\{ L_i \}},N)=Hom(j_{*}^{(I)}j^{(I)*}( \boxtimes
L_i),\Delta_{*}^{(I)}N);  $$ where $\boxtimes$ denotes the external
tensor product (it lives on $X^I$), all symbols like $f^{*}$ or $f_{*}$
denote the corresponding functors in the category of ${\cal
D}_{X^I}$-modules, and Hom is taken in that category as well. 

\begin{thm}([BD]). 

a)The sets from a) define a structure of a pseudo-tensor category on
${\cal M}(X)$. They are called $* I-operations$. 

b)The sets from b) define a structure of a pseudo-tensor category on
${\cal M}(X)$. They are called $chiral$ $operations$.The corresponding
category is denoted by ${\cal M}^{ch}$. 

\end{thm}

Since the usual operads are just pseudo-tensor categories with one
object,it is possible to speak about Lie algebras in pseudo-tensor
categories (they are pseudo-tensor subcategories which are functorial
images of the $Lie$-operad). 

\begin{dfn}. 

Lie algebras in the pseudo-tensor category a) are called
$Lie^{*}$-algebras.Lie algebras in the pseudo-tensor category b) are
called weak chiral algebras. 

\end{dfn}

If a weak chiral algebra contains a naturally defined unit (see
[BD],1.6.3) it is called $chiral$ $algebra$. 

The notion of a chiral algebra can be considered as a generalization of
the notion of a vertex algebra (see [FLM],[K],[Bo]) to the case of an
arbitrary curve $X$. An extensive treatment of chiral algebras from the
viewpoint of the theory of ${\cal D}$-modules can be found in [BD].

\subsection{ Remarks about q-deformed chiral algebras}

One can ask about q-deformed version of the notion of chiral algebra.
One cannot expect to have it for curves of the genuses higher than 1.
Currently there are exist few examples in genus zero case (see [FR]). 

In the traditional approach a vertex algebra (=chiral algebra on the
formal disk) is thought as a vector space $V$ equipped with a linear
map $V \otimes V \rightarrow V[[z,z^{-1}]]$ satisfying certain
properties.Equivalently, one can think of it as a lin ear map $V
\rightarrow (End V)[[z,z^{-1}]]$, $v \mapsto Y(v,z)$. One of the main
properties is locality: $(z-w)^N([Y(a,z),Y(b,w)])=0$ for a sufficiently
large integer $N$,and arbitrary $a,b \in V$.  One can try to replace
this condition by a more general one, like
$f(z/w)Y(a,z)Y(b,w)=Y(b,w)Y(a,z)$ where the matrix-valued function
$f(t)$ has singularities in the geometric series $q^n, n \in {\bf Z}$
or in a more general lattice. Thus the function $f$ becomes a new datum
of the theory. This idea was used i n [FR] where a preliminary
definition of a q-vertex algebra was given. 

One can try to interpret a q-deformed chiral algebra as a braided
version of a Lie algebra in a certain meromorphic pseudo-braided ${\bf
C^{*}}$-category. We hope to return to this topic in the future. 

\subsection{About G-vertex algebras}

R.Borcherds in [Bo] suggested a slightly different approach to vertex
algebras. The motivation is more or less as follows. Let us interpret
the vertex operator $Y(v,z)$ as $v^g$ where $g=exp(zL_{-1})$ and
$L_{-1}$ is the generator of the Virasoro algebra.  We can take as $g$
a more general element of the ``local Virasoro group'' $G$. Thus we can
make products like $v_1^{g_1},...,v_n^{g_n}$ corresponding to the
products of vertex operators. Applying this expression to $1 \in V$ we
obtain a space of operation s $P_{\delta_n}(G^n;V,...,V,V)$ where
$\delta_n$ is the only tree in ${\cal T}(n)$ which has no internal
edges. The space of operations can be informally interpreted as an
extension of $V[[g_1,...,g_n]]$ by $(g_i-g_j)^{-1}$. The same can be
done for any bracketing in $v_1...v_n$. These spaces are related by
the morphisms which roughly correspond to the restrictions to the
complements of the sub-divisors (one of the smaller diagonals in this
case). All this is compatible with the gluing of trees. Then one
obtains a $G-vertex$ $algebra$ which can be interpreted as an
associative algebra in the appropriate meromorphic pseudo-tensor
$G$-category. We refer to [Bo] for the precise definition and more
examples. In this way one can avoid $D$-modules and consider
multidimensional generalizations of vertex algebras related to various
groups $G$. 

We remark that a group $G$ as a datum can be useful even in the case
when a problem and the answer do no contain it. Let us return to the
example of conformal blocks in WZW model (see the end of Section 3.6).
Having a smooth curve $C$ of genus $0$ with $n $ marked points
$z_1,...,z_n$ and fixed (standard) local parameters at them one can
consider the corresponding affine Kac-Moody algebras $g_i^{\prime}$
``attached'' to these points.Let $g^{\prime}$ be an affine Kac-Moody
algebra corresponding to $z=0$ (the standard one),and $V_1,...,V_n$
highest weight representations having the same fixed central charge
$k$. Then under some mild conditions on $V_i$ and $k$ ,one can asign to
each plane tree $T \in {\cal T}(n)$ the ``fusion'' tensor product
$\odot_i V_i $ with the bracketing prescribed by $T$. It is known that
if $x$ belongs to the Virasoro algebra then $exp (x)$ acts on the
highest weight reprezentations. Thus we can make a tensor product
$\odot_i V_i(g_i)$ where $g_i$ are of the form $exp (z_iL_{-1})$. 
 In this way we get an analytic braided $G$-category where $G$ is a
1-parametric group generated by $L_{-1}$. Of course this analytic
braided category is equivalent to the usual braided category defined in
[KL]. In fact many proofs in [KL] use either $G$
 or the group $SL(2,{\bf C})$ with the Lie algebra generated by
$(L_{-1},L_0,L_1) \subset Vir$.

\vspace{15mm}

{\bf References}

\vspace{2mm}

[BD] A.Beilinson,V.Drinfeld, Chiral algebras,paper in preparation,1995. 

\vspace{2mm}

[BZ] J.Bernstein,A.Zelevinsky ,Representations of the group $GL(n,F)$
where $F$ is a non-archimedian local field.Russian
Math.Surveys,31:3(1976),1-68.  \vspace{2mm}

[Bo] R.Borcherds,Vertex algebras,preprint q-alg/9706008

\vspace{2mm}

[CP] V.Chari,A.Pressley, Quantum affine algebras and affine Hecke
algebras, Pacific J.Math.,174:2(1996),295-326. 

\vspace{2mm}

[De1] P.Deligne, Categories tannakiennes,Progress in Math.,vol.87
(1990),111-194. 

\vspace{2mm}

[De2] P.Deligne, Une description de categorie tressee,letter to
V.Drinfeld,1990.  \vspace{2mm}

[D1] V.Drinfeld, Quantum groups,Proceedings of ICM-86.American
Math.Soc.,1987. 

\vspace{2mm}

[D2] V.Drinfled,New realizations of Yangians and quantum affine
algebras,Russian Doklady,1985. 

\vspace{2mm}

[D3] V.Drinfeld, Degenerate affine Hecke algebras and Yangians, Funct.
anal. and appl. 20:1(1986),69-70. 

\vspace{2mm}

[FR] E.Frenkel,N.Reshetikhin ,Towards deformed chiral algebras,preprint
q-alg/9706023

\vspace{2mm}

[FLM] I.Frenkel,J.Lepowsky,A.Meurman,Vertex operator algebras and the
Monster,Pure and Appl.math.,vol.134,Academic Press,Boston,1988. 
\vspace{2mm}

[GK] V.Ginzburg,M.Kapranov Koszul duality for operads,Duke Math.J. 

\vspace{2mm}

[Gr] A.Grothendieck,Categories fibres and descente,Lecture Notes
Math.224,145-194

\vspace{2mm}

[K] V.Kac, Vertex algebras for beginners, University lecture series
vol.10,American Math.Soc.,1996. 

\vspace{2mm}

[Ka] M.Kapranov, The permutoassociahedron,Mac Lane's coherence theorem
and asymptotic zones for the KZ equation,J.Pure Appl.Alg. 85
(1993),119-142. 

\vspace{2mm}

[KL] D.Kazhdan,G.Lusztig ,Tensor structures arising from affine Lie
algebras,American J.Math.,1993-1994. 

\vspace{2mm}

[KM] M.Kontsevich,Yu.Manin ,Gromov-Witten classes,quantum cohomology
and enumerative geometry,Commun.Math.Phys.,164:3(1994),525-562. 

\vspace{2mm}

[KS] D.Kazhdan,Y.Soibelman Representations of affine quantum
algebras,Selecta Math.New Series,1996. 

\vspace{2mm} [L] G.Lusztig, Introduction to quantum
groups,Birkhauser,1993. 

\vspace{2mm}

[MS] G.Moore,N.Seiberg, Classical and quantum Conformal Field
Theory,Comm.Math.Phys.(1988). 

\vspace{2mm}

[S] F.Smirnov ,Remarks on deformed and undeformed KZ
equations,hep-th/9210051

\vspace{3mm}

$address$: Department of Mathematics,Kansas State
University,Manhattan,KS 66506,USA

$email$: soibel@math.ksu.edu

\end{document}